\documentclass[12pt, reqno]{amsart}

\usepackage{amssymb}  

\theoremstyle{plain}
\newtheorem{theorem}{Theorem}[section]
\newtheorem{lemma}[theorem]{Lemma}

\newtheorem{proposition}[theorem]{Proposition}
\newtheorem{corollary}[theorem]{Corollary}

\newtheorem{problem}[theorem]{Problem}

\theoremstyle{remark}
\newtheorem{example}[theorem]{Example}

\numberwithin{equation}{section}

\newcommand{\seclabel}[1]{\label{sec:#1}} 
\newcommand{\thmlabel}[1]{\label{thm:#1}} 
\newcommand{\lemlabel}[1]{\label{lem:#1}} 
\newcommand{\corlabel}[1]{\label{cor:#1}} 
\newcommand{\prplabel}[1]{\label{prp:#1}} 
\newcommand{\exmlabel}[1]{\label{exm:#1}} 
\newcommand{\eqnlabel}[1]{\label{eqn:#1}} 

\newcommand{\secref}[1]{\ref{sec:#1}} 
\newcommand{\thmref}[1]{\ref{thm:#1}} 
\newcommand{\lemref}[1]{\ref{lem:#1}} 
\newcommand{\corref}[1]{\ref{cor:#1}} 
\newcommand{\prpref}[1]{\ref{prp:#1}} 
\newcommand{\exmref}[1]{\ref{exm:#1}} 
\newcommand{\eqnref}[1]{\ref{eqn:#1}} 

\newcommand{\peqref}[1]{(\eqnref{#1})} 

\newcommand{\Aut}{\mathrm{Aut}}     
\newcommand{\Mlt}{\mathrm{Mlt}}     
\newcommand{\Fix}{\mathrm{Fix}}     
\newcommand{\Ker}{\mathrm{Ker}}         
\newcommand{\LNuc}{N_{\lambda}}     
\newcommand{\MNuc}{N_{\mu}}                 
\newcommand{\RNuc}{N_{\rho}}        

\newcommand{\setof}[2]{\{#1\,|\,#2\}}   

\newcommand{\sbl}[1]{\langle#1\rangle}          
\newcommand{\iv}{^{-1}}             

\title[When is the commutant of a Bol loop a subloop?]
{When is the commutant of a Bol loop a subloop?}

\author[M.~K.~Kinyon]{Michael~K.~Kinyon}
\address{Department of Mathematical Sciences \\
Indiana University South Bend \\
South Bend, IN 46634 USA}
\email{mkinyon@iusb.edu}
\urladdr{http://mypage.iusb.edu/\symbol{126}mkinyon}
\author[J.~D.~Phillips]{J.~D.~Phillips}
\address{Department of Mathematics \& Computer Science \\
Wabash College \\
Craw\-fords\-ville, IN 47933 USA}
\email{phillipj@wabash.edu}
\urladdr{http://www.wabash.edu/depart/math/faculty.html{\#}Phillips}
\author[P. Vojt\v{e}chovsk\'{y}]{Petr Vojt\v{e}chovsk\'{y}}
\address{Department of Mathematics \\
University of Denver \\
2360 S Gaylord St \\
Denver, CO 80208 USA}
\email{petr@math.du.edu}
\urladdr{http://www.math.du.edu/\symbol{126}petr}

\date{\today}


\begin{document}

\begin{abstract}
A left Bol loop is a loop satisfying $x(y(xz)) = (x(yx))z$. The commutant of a
loop is the set of elements which commute with all elements of the loop.
In a finite Bol loop of odd order or of order $2k$, $k$ odd, the commutant is
a subloop. We investigate conditions under which the commutant of a Bol loop is not a
subloop. In a finite Bol loop of order relatively prime to $3$, the commutant generates
an abelian group of order dividing the order of the loop. This generalizes a well-known
result for Moufang loops.
After describing all extensions of a loop $K$ such that $K$ is in the
left and middle nuclei of the resulting loop, we show how to construct classes
of Bol loops with non-subloop commutant. In particular, we obtain all Bol loops
of order $16$ with non-subloop commutant.
\end{abstract}

\maketitle

\section{Introduction}
\seclabel{intro}

A \emph{loop} $(Q,\cdot)$ is a set $Q$ with a binary operation $\cdot$ such
that there is a neutral element $1\in Q$ satisfying $1\cdot x = x\cdot 1 = x$
for all $x\in Q$, and such that for each $a$, $b\in Q$ the equations $a\cdot x
= b$ and $y\cdot a = b$ have unique solutions $x$, $y\in Q$. We write $xy$
instead of $x\cdot y$, and reserve $\cdot$ to have lower priority than
juxtaposition among factors to be multiplied---for instance, $x\cdot yz$
stands for $x(yz)$.

The \emph{commutant} (also known as the \emph{centrum}, \emph{Moufang center}
or \emph{semicenter}) of a loop $Q$ is the set
\begin{displaymath}
C(Q) = \setof{c\in Q}{cx = xc\text{ for every }x\in Q}.
\end{displaymath}
In a group, or even a Moufang loop, the commutant is a subloop, but this does
not need to be the case in general. When $Q$ is a loop and $C(Q)$ is not a
subloop of $Q$, we say that $Q$ \emph{has a non-subloop commutant}.

Given an element $a$ of a loop $Q$, we denote by $L_a$ the left translation of
$Q$ by $a$, \textit{i.e.}, $bL_a = ab$. Similarly, $R_a$ is the right translation
$bR_a = ba$. The commutant is obviously characterized as
$C(Q) = \setof{c\in Q}{L_c = R_c}$.
The permutation group $\sbl{L_a,\,R_a\,|\,a\in Q}$ is known as the
\emph{multiplication group} of $Q$, and will be denoted by $\Mlt(Q)$.
We also use the notations $R : Q\to \Mlt(Q) ; a\mapsto R_a$ and
$L : Q\to \Mlt(Q); a\mapsto L_a$.

A loop is \emph{left Bol} if it satisfies
\begin{displaymath}
    x(y \cdot xz) = (x \cdot yx)z  \tag{\textsc{Bol}}
\end{displaymath}
for all $x,y,z$, or equivalently, if $L_x L_y L_x = L_{x\cdot yx}$ for all $x,y$.
Right Bol loops are defined by the mirror of (\textsc{Bol}). We will consider
only left Bol loops in this paper, and henceforth refer to them simply as Bol loops.
Note that much of the literature on Bol loops (\textit{e.g.}, \cite{DR-loop}) considers
right Bol loops, and hence results need to be translated appropriately.

The main purposes of this paper are to introduce general constructions of Bol loops
with non-subloop commutants, and to shed light on the structure of the commutant of
a Bol loop, of the subloop generated by the commutant, and of the subgroup
$\sbl{L_c\,|\,c\in C(Q)}$ of $\Mlt(Q)$.

The commutant of a finite Bol loop of odd order is a subloop \cite{KP}. In
Theorem \thmref{order4k}, we show that the commutant of a finite Bol loop
of order $2k$, $k$ odd, is a subloop. Thus finite Bol loops with non-subloop
commutant have order divisible by $4$ (Corollary \corref{order4k}). When the
commutant of a Bol loop is a subloop, it is a commutative Moufang loop.
In case $Q$ is a Bol loop with a non-subloop commutant, it is natural to
consider the structure of the subloop $\sbl{C(Q)}$. In Corollary \corref{ab-grp2},
we show that if $Q$ is finite of order relatively prime to $3$, then
$\sbl{C(Q)}$ is an abelian group of order dividing $|Q|$. This generalizes a well-known
result about commutants of Moufang loops. We then turn to constructions. In \S\secref{leftnuc},
we describe all extensions $Q$ of a group $K$ such that $K$ is contained in the
left and middle nuclei of the resulting loop. In the next two sections, we
specialize this to construct examples of Bol loops with non-subloop commutants.
In \S\secref{semi}, we consider the special case of a semidirect product
(split extension), and in \S\secref{trivact}, we consider the special case
where the action of $Q/K$ on $K$ is trivial. When we restrict to low orders,
the two constructions give $20$ of the $21$ known Bol loops of order less than or
equal to $16$ with non-subloop commutant. We finish with another construction
which yields the remaining such loop.

We conclude this introduction with a review of some basic facts regarding loops
in general and Bol loops in particular. The standard references \cite{Br, Pf}
provide adequate general background in loop theory, while the latter reference,
\cite[Chap. VI]{Kiechle}, and \cite{DR-loop} give specific details regarding Bol loops.

For a loop $Q$, the following subgroups (associative subloops) are of interest:
\begin{center}
\begin{tabular}{rlrcl}
$\bullet$ & the \emph{left nucleus} &
$\LNuc(Q)$ & $=$ & $\setof{a\in Q}{a\cdot xy = ax\cdot y, \forall x,y\in Q}$ \\
$\bullet$ & the \emph{middle nucleus} &
$\MNuc(Q)$ & $=$ & $\setof{a\in Q}{x\cdot ay = xa\cdot y, \forall x,y\in Q}$ \\
$\bullet$ & the \emph{right nucleus} &
$\RNuc(Q)$ & $=$ & $\setof{a\in Q}{x\cdot ya = xy\cdot a, \forall x,y\in Q}$ \\
$\bullet$ & the \emph{nucleus} &
$N(Q)$ & $=$ & $\LNuc(Q)\cap \MNuc(Q)\cap \RNuc(Q)$ \\
$\bullet$ & the \emph{center} &
$Z(Q)$ & $=$ & $N(Q)\cap C(Q)$ \\
\end{tabular}
\end{center}
The center is a normal subloop. In a Bol loop $Q$, the left and middle nuclei
coincide, $\LNuc(Q) = \MNuc(Q)$, and we shall just refer to this as the left nucleus.
The left nucleus of a Bol loop is a normal subloop, but does not necessarily
coincide with the right nucleus, nor is the right nucleus necessarily normal.
In addition, $Z(Q) = C(Q)\cap \LNuc(Q)$. Indeed, for $c\in C(Q)\cap \LNuc(Q)$
and $x,y\in Q$, $xy\cdot c = c\cdot xy = cx\cdot y = xc\cdot y = x\cdot cy
= x\cdot yc$. Thus $c\in \RNuc(Q)$ and so $c\in Z(Q)$.

Let $Q$ be a Bol loop. Then $Q$ is \emph{power-associative}, that is, for each
$x\in Q$, the subloop $\sbl{x}$ generated by $x$ is a subgroup. In particular,
there exists $x\iv\in Q$ such that $xx\iv = x\iv x = 1$. In addition, $Q$ is
\emph{left power alternative}, that is,
\[
x^m \cdot x^n y = x^{m+n} y  \tag{\textsc{LPA}}
\]
for all integers $m,n$. Equivalently, $L_x^m = L_{x^m}$ for every  integer $m$.
An element $a\in Q$ is said to be \emph{right power alternative} if
$R_a^m = R_{a^m}$ for every integer $m$. Not every element of a nonMoufang
Bol loop $Q$ is right power alternative.

For a Bol loop $Q$, a subset $S$ is a subloop if and only if it is closed under
both multiplication and inversion. Indeed, if $a,b\in S$, then the unique solutions
of the equations $ax=b$ and $ya=b$ are $x=a\iv b\in S$ and
$y = a\iv (ab\cdot a\iv)\in S$, respectively.

\begin{lemma}
\lemlabel{order}
Let $Q$ be a finite loop, and let $S\subseteq Q$ be a subloop.
\begin{enumerate}
\item[(i)] If $R |_S : S\to \Mlt(Q)$ is a homomorphism, then $|S|$ divides $|Q|$.
\item[(ii)] If $L |_S : S\to \Mlt(Q)$ is a homomorphism, then $|S|$ divides $|Q|$.
\end{enumerate}
\end{lemma}

\begin{proof}
Fix $x,y\in Q$, and suppose $xS\cap yS \neq \emptyset$. Then there exist
$s_1, s_2\in S$ such that $xs_1 = ys_2$, and so $x = y R_{s_2} R_{s_1}\iv
= y R_{s_3} = y s_3$ where $s_3 s_1 = s_2$. Thus for all $s\in S$,
$xs = y R_{s_3}R_s = y \cdot s_3 s \in yS$, and so $xS \subseteq yS$.
The other inclusion follows similarly, and so $xS = yS$. Therefore
$\setof{xS}{x\in Q}$ is a partition of $Q$. This establishes (i), and
the proof of (ii) is similar, using left cosets instead of right cosets.
\end{proof}

\begin{corollary}
\corlabel{order}
Let $Q$ be a finite loop.
\begin{enumerate}
\item[(i)] If $S$ is a subgroup of $\RNuc(Q)$, then $|S|$ divides $|Q|$.
\item[(ii)] If $S$ is a subgroup of $\LNuc(Q)$, then $|S|$ divides $|Q|$.
\end{enumerate}
If, in addition, $Q$ is a Bol loop, then for each $x\in Q$,
$|\sbl{x}|$ divides $|Q|$.
\end{corollary}

\begin{proof}
Parts (i) and (ii) follow from the respective parts of Lemma \lemref{order}.
The remaining assertion is well-known \cite{Pf,DR-loop}, and follows from (\textsc{LPA}).
\end{proof}

It is not known if the order of an arbitrary subloop of a finite Bol loop
divides the order of $Q$.

\section{Structure of the Commutant}
\seclabel{structure}

\emph{Throughout this section}, let $Q$ denote a Bol loop with commutant $C(Q)$.

\begin{lemma}
\lemlabel{ab}
If $a\in C(Q)$, then:
\begin{enumerate}
\item[(i)]\quad $\sbl{a}\subset C(Q)$,
\item[(ii)]\quad $a$ is right power alternative.
\end{enumerate}
In addition, if $b\in C(Q)$, then:
\begin{enumerate}
\item[(iii)]\quad $a^2 b\in C(Q)$,
\item[(iv)]\quad $a^k b^{\ell}\cdot a^m b^n = a^{k+m} b^{\ell + n}$ for all integers $k,\ell,m,n$.
\end{enumerate}
\end{lemma}

\begin{proof}
For (i) and (iii), see Lemmas 2.2 and 2.1 of \cite{KP}. Part (ii)
follows from (i) and (\textsc{LPA}). Let us prove (iv):
\begin{align*}
a^k b^{\ell}\cdot a^m b^n &= a^{-m}\cdot a^m (a^k b^{\ell}\cdot a^m b^n)
= a^{-m} [ a^m (a^k b^{\ell}\cdot a^m)\cdot b^n] \\
&= a^{-m}\cdot (a^{2m}\cdot a^k b^{\ell})b^n = a^{-m} (a^{k+2m} b^{\ell}\cdot b^n) \\
&= a^{-m}\cdot a^{k+2m} b^{\ell + n} = a^{k+m} b^{\ell + n}
\end{align*}
using (\textsc{LPA}), (\textsc{Bol}), (i), (\textsc{LPA}) twice, (ii), and
(\textsc{LPA}) once more.
\end{proof}

Part (i) of the lemma immediately implies the following.

\begin{corollary}
\corlabel{sub-mult}
$C(Q)$ is a subloop if and only if it is closed under multiplication.
\end{corollary}

We introduce here some useful notation. For an integer $n > 1$, let
\[
C_n(Q) := \setof{c\in C(Q)}{c\text{ has finite order relatively prime to }n} .
\]
Obviously, if $m>1$ divides $n$, $C_n(Q)\subseteq C_m(Q)$.

\begin{lemma}
\lemlabel{cn-ab}
Let $n > 1$ be an integer.
\begin{enumerate}
\item[(i)]\quad If $a\in C_n(Q)$, then $\sbl{a}\subseteq C_n(Q)$.
\item[(ii)]\quad If $a,b\in C_n(Q)$, then $ab$ has order relatively prime to $n$.
\item[(iii)]\quad If $a,b\in C_n(Q)$, then $a^2 b\in C_n(Q)$.
\item[(iv)]\quad If $m > 1$ divides $n$, and if $C_m(Q)$ is a subloop, then $C_n(Q)$ is a subloop.
\end{enumerate}
\end{lemma}

\begin{proof}
Part (i) is immediate from Lemma \lemref{ab}(i), part (ii) follows from
Lemma \lemref{ab}(iv), and (iii) follows from (ii) and Lemma \lemref{ab}(iii).
Finally, (iv) follows from (ii) and $C_n(Q)\subseteq C_m(Q)$.
\end{proof}

The following is a mild generalization of the main result of \cite{KP}.

\begin{corollary}
\corlabel{odd}
For each integer $m > 0$, $C_{2m}(Q)$ is a subloop of $Q$.
\end{corollary}

\begin{proof}
By Lemma \lemref{cn-ab}(iv), it is enough to show that $C_2(Q)$
is a subloop. By Lemma \lemref{cn-ab}(i), $C_2(Q)$ is closed under inversion.
If $a,b\in C_2(Q)$, then by Lemma \lemref{cn-ab}(i), there exists $c\in C_2(Q)$
such that $c^2 = a$, and so $ab = c^2 b\in C_2(Q)$, using Lemma \lemref{cn-ab}(iii).
\end{proof}

\begin{lemma}
\lemlabel{rightnuc}
For $c\in C(Q)$, $c^2 \in \LNuc(Q)$ if and only if $c\in \RNuc(Q)$.
\end{lemma}

\begin{proof}
For $x$, $y\in Q$, we compute
\begin{equation}
\eqnlabel{rn1}
c^2\cdot xy = c(c\cdot xy) = (xy\cdot c)c
\end{equation}
and
\begin{equation}
\eqnlabel{rn2}
c^2x\cdot y = (c\cdot xc)y = c(x\cdot cy) = (x\cdot yc)c .
\end{equation}
Now $c^2\in \LNuc(Q)$ iff the left hand sides of
\peqref{rn1} and \peqref{rn2} are equal, while $c\in \RNuc(Q)$
iff the right hand sides are equal. The result follows.
\end{proof}

\begin{corollary}
\corlabel{oddcenter}
If $C(Q) = C_2(Q)$, then $Z(Q) = C(Q)\cap \RNuc(Q)$.
\end{corollary}

\begin{proof}
In this case, the lemma gives $C(Q)\cap \RNuc(Q) = C(Q)\cap \LNuc(Q) = Z(Q)$.
\end{proof}

\begin{corollary}
\corlabel{order2}
If $S\subseteq \setof{c\in C(Q)}{c^2 \in \LNuc(Q)}$, then $S$
generates an abelian subgroup of $\RNuc(Q)$.
\end{corollary}

\begin{theorem}
\thmlabel{order4k}
Let $Q$ be a finite Bol loop of order $2k$ where $k$ is odd. Then $C(Q)$
is a subloop of $Q$.
\end{theorem}

\begin{proof}
If $C(Q) = C_2(Q)$, then the result follows from
Corollary \corref{odd}. Thus assume $1 \neq a\in C(Q)$ has order $2$.
Fix $b\in C(Q)$ of even order $2m$ where $m$ divides $k$. Then $b^m$ has
order $2$. If $b^m \neq a$, then by Corollary \corref{order2},
$\sbl{a,b^m}$ is an abelian subgroup of $\RNuc(Q)$ of order $4$. But then
by Corollary \corref{order}(i), $4$ divides $|Q|$, a contradiction. Thus $b^m = a$.
Set $c = b^{m+1}$ and note that $b = ac$. Hence $c^m = a^m b^m = aa = 1$.

Summarizing, every element of $C(Q)$ can be written uniquely in the form
$a^i c$ where $i\in \{0,1\}$ and where $c\in C_2(Q)$. For
$c_1, c_2\in C_2(Q)$, $c_1 c_2\in C_2(Q)$
(Corollary \corref{odd}), and
$a^i c_1\cdot a^j c_2 = a^{i+j}\cdot c_1 c_2$, $i,j\in \{0,1\}$, by
Lemma \lemref{ab}(iv). Thus $C(Q)$ is closed under multiplication, and
so by Corollary \corref{sub-mult}, it is a subloop.
\end{proof}

\begin{corollary}
\corlabel{order4k}
If $Q$ is a finite Bol loop with non-subloop commutant, then $4$ divides $|Q|$.
\end{corollary}

In Proposition \prpref{order4n}, we will show that for each
integer $n>2$, there exists a Bol loop of order $4n$ with
non-subloop commutant.

\section{Commutant elements of order prime to $3$}
\seclabel{prime3}

We now proceed to show that in a Bol loop, commutant elements of order
relatively prime to $3$ generate an abelian group. This generalizes
the well-known result that in a Moufang loop, commutant elements of order
relatively prime to $3$ lie in the center.
There are nonassociative commutative Moufang loops, the smallest
being of order $81$, and so the assumptions on the orders of elements
or of loops are necessary.

We adopt the following convention: for elements $a_1, a_2, \ldots, a_n$
of a loop $Q$, $a_1 a_2\cdots a_n$ will denote the left-associated product
$(\cdots (a_1 a_2)\cdots ) a_n$.

\begin{theorem}
\thmlabel{righttrans}
Let $Q$ be a Bol loop, let $A\subseteq C(Q)$, and suppose that for each
$a,b\in A$, $R_a R_b = R_{ab}$. Then the subloop
$\sbl{A}$ is an abelian subgroup of $Q$, and
$R |_{\sbl{A}} : \sbl{A}\to \Mlt(Q)$ is a homomorphism.

If, in addition, $Q$ is finite, then $|\sbl{A}|$ divides $|Q|$.
\end{theorem}

\begin{proof}
Since $R_a R_b = R_b R_a$ for all $a,b\in A$, we may freely rearrange
products of right translations from $A$. For $n>0$, let $a_1,\ldots, a_n\in A$.
We will verify
\begin{equation}
\eqnlabel{newkey}
R_{a_1 a_2\cdots a_n} = R_{a_1} R_{a_2} \cdots R_{a_n}
\end{equation}
by induction on $n$. By hypothesis, \peqref{newkey} holds for $1\leq n\leq 2$.
Suppose $n > 2$ and that \peqref{newkey} holds for $n-1$. Then
\begin{align*}
a_{n-1}\cdot \lbrack x R_{a_1} \cdots R_{a_{n-2}} R_{a_{n-1}} R_{a_n} \rbrack
&= a_{n-1} \lbrack x R_{a_1} \cdots R_{a_{n-2}} \cdot a_{n-1} a_n \rbrack \\
&= \left( a_{n-1} \cdot x R_{a_1} \cdots R_{a_{n-2}} R_{a_{n-1}}\right) a_n \\
&= \left( a_{n-1} \cdot x R_{a_1 \cdots a_{n-2} a_{n-1}}\right) a_n \\
&= \left( a_{n-1} \cdot x \lbrack a_{n-1}\cdot (a_1 \cdots a_{n-2})\rbrack \right) a_n \\
&= \lbrack ( a_{n-1} \cdot x a_{n-1} ) \cdot (a_1 \cdots a_{n-2}) \rbrack a_n \\
&= ( a_{n-1} \cdot x a_{n-1} ) (a_1 \cdots a_{n-2} a_n ) \\
&= a_{n-1} \cdot x \lbrack a_{n-1} \cdot (a_1 \cdots a_{n-2} a_n )\rbrack \\
&= a_{n-1} \cdot x (a_1 \cdots a_{n-2} a_{n-1} a_n ) . \\
\end{align*}
Here we are using, in succession, $R_{a_{n-1}} R_{a_n} = R_{a_{n-1} a_n}$,
(\textsc{Bol}), the induction hypothesis, $a_{n-1}\in C(Q)$, (\textsc{Bol})
again, the induction hypothesis again, (\textsc{Bol}) once more,  $a_{n-1}\in C(Q)$
again, and $R_{a_n} R_{a_{n-1}} = R_{a_{n-1}} R_{a_n}$. Cancelling $a_{n-1}$, we
obtain \peqref{newkey} for $n$.

Since $a_1 a_2 \cdots a_n = 1 R_{a_1} R_{a_2}\cdots R_{a_n}$, and since
$R_{a_i} R_{a_j} = R_{a_j} R_{a_i}$ for all $i,j$, it follows that
the expression $a_1 a_2 \cdots a_n$ is invariant under all reassociations
and rearrangements. Thus $\sbl{A}$ is an abelian group, and the homomorphism
assertion follows from \peqref{newkey}.

The remaining claim follows from Lemma \lemref{order}(i).
\end{proof}

\begin{lemma}
\lemlabel{rearr}
Let $Q$ be a Bol loop, and suppose $a,b\in C(Q)$ satisfy $R_a R_b = R_b R_a$.
Then $R_a R_b = R_{ab}$.
\end{lemma}

\begin{proof}
We compute
\[
a(xa\cdot b) = (xa\cdot b)a = (xa\cdot a)b = (a\cdot xa)b = a(x\cdot ab) ,
\]
using (\textsc{Bol}) in the last step. Cancelling $a$, we have the desired result.
\end{proof}

\begin{lemma}
\lemlabel{cubes}
Let $Q$ be a Bol loop, and let $a$, $b$, $c\in C(Q)$. Then:
\begin{enumerate}
\item[(i)]\quad for all $x\in Q$, $xb\cdot a^3 = x a^3 \cdot b = x\cdot a^3 b$,
\item[(ii)]\quad for all $x\in Q$, $x^3 a\cdot b = x^3 b\cdot a = x^3\cdot ab$.
\end{enumerate}
\end{lemma}

\begin{proof}
By Lemma \lemref{ab}, $a(b\cdot ax) = (a\cdot ba)x = a^2 b\cdot x = x \cdot
a^2 b$. Thus $a^3 (ax\cdot b) = a^3 (b\cdot ax) = a^2\cdot a(b\cdot ax) =
a^2(x\cdot a^2 b) = (a^2\cdot xa^2)b = (a^3\cdot ax) b$. Replacing $x$ with
$a\iv x$ and using $a^3 , b\in C(Q)$, we have the first equality of (i).
The second equality follows from Lemma \lemref{rearr}.

Next, we compute $a\cdot x^2 a = a^3\cdot a^{-2}(x\cdot xa) = a^3(a\iv\cdot
(x\cdot ax)a\iv)$, and so $a(x^2\cdot ab) = (a\cdot x^2 a)b = a^3(a\iv\cdot
(x\cdot ax)a\iv)\cdot b = a^3\cdot (a\iv\cdot (x\cdot ax)a\iv)b = a^3\cdot a\iv
((x\cdot ax)\cdot a\iv b) = a\cdot a(x\cdot a(x\cdot a\iv b))$, using (i) in
the third equality. Cancelling $a$ on the left and multiplying by $x$ on the left, we have
\begin{equation}\eqnlabel{Aux1}
    x^3\cdot ab = x\cdot a(x\cdot a(x\cdot a\iv b))
    = (x\cdot (a\cdot xa)x)\cdot a\iv b.
\end{equation}
Since $x\cdot (a\cdot xa)x = x\cdot a(x\cdot ax) = x\cdot a(x^2 a) = x\cdot
x^2a^2 = x^3a^2 = a \cdot x^3 a$, we can rewrite (\eqnref{Aux1}) as $x^3\cdot
ab = (a \cdot x^3 a)\cdot a\iv b = a\cdot x^3 b = x^3 b\cdot a$, and (ii)
follows.
\end{proof}

\begin{corollary}
\corlabel{c3-ab-grp}
Let $Q$ be a Bol loop. For each positive integer $n$, $\sbl{C_{3n}(Q)}$ is an abelian group,
and $R |_{\sbl{C_{3n}(Q)}} : \sbl{C_{3n}(Q)}\to \Mlt(Q)$ is a homomorphism.
If, in addition, $Q$ is finite, then $|\sbl{C_{3n}(Q)}|$ divides $|Q|$.
\end{corollary}

\begin{proof}
By Lemma \lemref{cubes}(i), $R_a R_b = R_{ab}$ for all $a,b\in C_{3n}(Q)$, and so
Theorem \thmref{righttrans} applies with $A = C_{3n}(Q)$.
\end{proof}

\begin{corollary}
\corlabel{ab-grp}
Let $Q$ be a Bol loop such that $C(Q) = C_{3n}(Q)$ for some integer $n>0$.
Then $\sbl{C(Q)}$ is an abelian group, and
$R |_{\sbl{C(Q)}} : \sbl{C(Q)}\to \Mlt(Q)$ is a homomorphism.
If, in addition, $Q$ is finite, then $|\sbl{C(Q)}|$ divides $|Q|$.
\end{corollary}

\begin{corollary}
\corlabel{ab-grp2}
Let $Q$ be a finite Bol loop of order relatively prime to $3$.
Then $\sbl{C(Q)}$ is an abelian group,
$R |_{\sbl{C(Q)}} : \sbl{C(Q)}\to \Mlt(Q)$ is a homomorphism,
and $|\sbl{C(Q)}|$ divides $|Q|$.
\end{corollary}

Note that one cannot replace right translations with left
translations in Corollary \corref{ab-grp}, for otherwise
$C(Q)$ would necessarily be a subloop.

Recall that any two elements of a Moufang loop generate a group, \textit{i.e.},
Moufang loops are \emph{diassociative}. It is well known that nonMoufang Bol loops
are not diassociative. However, after seeing the calculations in the proof of Lemma
\lemref{cubes}, the reader might wonder if in a Bol loop, two elements
generate a group if one of the two elements is in the commutant. The answer is
``no'':

\begin{example} Let $Q$ be the Bol loop
\begin{displaymath}
\begin{array}{c|cccccccc}
    &1&2&3&4&5&6&7&8\\
    \hline
1& 1& 2& 3& 4& 5& 6& 7& 8\\
2& 2& 1& 4& 3& 6& 5& 8& 7\\
3& 3& 4& 1& 2& 7& 8& 5& 6\\
4& 4& 3& 2& 1& 8& 7& 6& 5\\
5& 5& 6& 7& 8& 1& 2& 3& 4\\
6& 6& 5& 8& 7& 2& 1& 4& 3\\
7& 7& 8& 5& 6& 4& 3& 2& 1\\
8& 8& 7& 6& 5& 3& 4& 1& 2
\end{array}
\end{displaymath}
Then $Z(Q) = \LNuc(Q) = \{1,2\}$, $C(Q) = \RNuc(Q) = \{1,2,3,4\}$, and $\sbl{4,5} = Q$.
\end{example}

We do not know the answer to the following.

\begin{problem}
Does there exist a finite Bol loop of order relatively prime to $3$ such that
the commutant is not contained in the right nucleus?
\end{problem}

\section{Left Nuclear Extensions of Bol Loops}
\seclabel{leftnuc}

Let $(Q,\cdot)$, $(K,\cdot)$, $(E,*)$ be loops. Then $Q$ is an \emph{extension
of $K$ by $E$} if $K$ is a normal subloop of $Q$ and $Q/K=E$. We can then
identify $E$ with a subset of $Q$---in fact, with a transversal of $Q/K$---and
assume without loss of generality that $1_E=1_Q=1$. Given $a$, $b\in E$, there
is then a unique $f(a,b)\in K$ such that $ab = f(a,b)(a*b)$. The map $f:E\times
E\to K$ thus obtained satisfies $f(a,1)=f(1,a)=1$. We will call a map with the
property $f(a,1)=f(1,a)=1$ a \emph{cocycle}.

\begin{theorem}
\thmlabel{AllExts}
Let $(K,\cdot)$, $(E,*)$ be loops, and $Q$ an extension of $K$ by $E$ such that
$K\le \LNuc(Q)\cap \MNuc(Q)$. Then $K$ is a group, there is a map
$f:E\times E\to K$ satisfying $f(1,a)=f(a,1)=1$, and a map
$\tau:E\to \Aut(K)$ satisfying $\tau_1=1$, such that $Q$ is isomorphic
to $K\times E$ with multiplication
\begin{equation}\eqnlabel{AllExts}
    (u,a)(v,b) = (u\tau_a(v)f(a,b),a*b)
\end{equation}
for every $a$, $b\in E$, $u$, $v\in K$.

Conversely, given a group $(K,\cdot)$, a loop $(E,*)$, a cocycle $f:E\times
E\to K$, and a map $\tau:E \to \Aut(K)$ with $\tau_1=1$, the loop
$Q=K\times E$ with multiplication \emph{\peqref{AllExts}} is an extension of $K$
by $E$, and $K\le \LNuc(Q)\cap \MNuc(Q)$.
\end{theorem}

\begin{proof}
$K$ is obviously a group since it is a subloop of two nuclei. Let $f:K\times K\to E$
be the cocycle described above. For $a\in E$, let $\tau_a : K\to K$ by defined by
$\tau_a(u) a = au$ for every $u\in K$. Since $\tau_a$ is the restriction of the
inner mapping $L_aR_a^{-1}$ of $Q$ to $K$, and since $K$ is normal in $Q$,
$\tau_a$ is a bijection of $K$.

We claim that $\tau_a$ is a homomorphism. We have
$\tau_a(uv)=\tau_a(u)\tau_a(v)$ if and only if $\tau_a(uv)a =
\tau_a(u)\tau_a(v)\cdot a = \tau_a(u)\cdot \tau_a(v)a$ (since $K\le
\LNuc(Q)$) if and only if $a(uv)=\tau_a(u)\cdot av$ if and only if
$a\cdot uv = \tau_a(u)a\cdot v= au\cdot v$ (again by $K\le \LNuc(Q)$). But
$a\cdot uv = au\cdot v$, since $u\in \MNuc(Q)$. Thus $\tau_a$ is a homomorphism, and
$\tau:a\mapsto\tau_a$ is a map $E \to \Aut(K)$. We see right away that
$\tau_1=1$.

Every element of $Q$ can be expressed uniquely as $ua$, where $u\in K$ and
$a\in E$ (since $E$ is a transversal of $Q/K$). For $u$, $v\in K$, $a$, $b\in E$,
we have: $ua\cdot vb = u(a\cdot vb)$ (since $u\in K\le \LNuc(Q)$),
$u(a\cdot vb) = u(av\cdot b)$ (since $v\in K\le \MNuc(Q)$). As $\tau_a(v)\in K$
and $f(a,b)\in K$, we have $u(av\cdot b) = u(\tau_a(v)a\cdot b) =
u(\tau_a(v)\cdot ab) = u\tau_a(v)\cdot ab = u\tau_a(v)\cdot f(a,b)(a*b) =
u\tau_a(v)f(a,b)\cdot (a*b)$.

For the converse, it is easy to see that $Q=K\times E$ with multiplication
(\eqnref{AllExts}) is a loop, $K\unlhd Q$, and $Q/K=E$. We have
\begin{eqnarray*}
(u,1)(v,b)\cdot (w,c) &=& (uv\tau_b(w)f(b,c),bc) = (u,1)\cdot (v,b)(w,c), \text{ and}\\
(v,b)(u,1)\cdot (w,c) &=& (v\tau_b(u)\tau_b(w)f(b,c),bc) = (v\tau_b(uw)f(b,c),bc) \\
&=& (v,b)\cdot (u,1)(w,c) ,
\end{eqnarray*}
and hence $K\le \LNuc(Q)\cap \MNuc(Q)$.
\end{proof}

Since the left and middle nuclei coincide in Bol loops, we have the following.

\begin{corollary}
Let $K$ be a group, $E$ a Bol loop. Assume that $Q$ is a Bol loop which is an
extension of $K$ by $E$, and that $K\le \LNuc(Q)$. Then
the multiplication in $Q$ is given by \emph{\peqref{AllExts}} for some cocycle $f$
and a map $\tau:E\to\Aut(K)$ satisfying $\tau_1=1$.
\end{corollary}

Denote the extension of $K$ by $E$ constructed as in \peqref{AllExts} by
$Q = Q(K,E,\tau,f)$. We are now going to give conditions on $f$ and $\tau$ that make
$Q$ into a Bol loop.

\begin{theorem}
\thmlabel{LB}
Let $K$ be a group, $E$ a Bol loop, $f:E\times E\to K$ a
cocycle and $\tau:E\to\Aut(K)$ a map satisfying $\tau_1=1$,
and set $Q=Q(K,E,\tau,f)$. Then:
\begin{enumerate}
\item[(i)] $Q$ is a Bol loop if and only if
\begin{eqnarray}
\eqnlabel{LBExtension1}
\tau_a(f(b,a))f(a,ba)f(a\cdot ba,c) &=&
\tau_a \tau_b(f(a,c))\tau_a(f(b,ac))f(a,b\cdot ac), \\
\eqnlabel{LBExtension2}
\tau_a(f(b,a))f(a,ba)\tau_{a\cdot ba}(w) &=&
\tau_a \tau_b \tau_a(w) \tau_a(f(b,a))f(a,ba)
\end{eqnarray}
for every $w\in K$ and $a, b, c\in E$.
\item[(ii)] $(w,c)\in Q$ belongs to $\RNuc(Q)$ if and only if
$c\in \RNuc(E)$ and
\begin{equation}
\eqnlabel{RNExtension}
f(a,b) \tau_{ab}(w) f(ab,c) = \tau_a \tau_b(w) \tau_a(f(b,c)) f(a,bc)
\end{equation}
for all $a,b\in E$.
\item[(iii)] $Q$ is a group if and only if $E$ is a group and
\begin{eqnarray}
\eqnlabel{GExtension1}
\tau_a (f(b,c))f(a,bc) &=& f(a,b) f(ab,c), \\
\eqnlabel{GExtension2}
\tau_a \tau_b(w) f(a,b) &=& f(a,b)\tau_{ab}(w)
\end{eqnarray}
for every $w\in K$ and $a, b, c\in E$.
\item[(iv)] $(u,a)\in Q$ belongs to $C(Q)$ if and only if $a\in C(E)$ and
\begin{eqnarray}
\eqnlabel{InC1}
\tau_a(v) &=& u\iv vu \\
\eqnlabel{InC2}
\tau_b(u) &=& u f(a,b) f(b,a)\iv
\end{eqnarray}
for every $v\in K$, $b\in E$.
\end{enumerate}
\end{theorem}

\begin{proof}
For (i): By straightforward calculation with $x=(u,a)$, $y=(v,b)$,
$z=(w,c)$ substituted into the Bol identity $(x\cdot yx)z = x(y\cdot xz)$,
we obtain that $Q$ is a Bol loop if and only if
\begin{equation}
\eqnlabel{LBExt-tmp}
\tau_a(f(b,a))f(a,ba)\tau_{a\cdot ba}(w)f(a\cdot ba,c)
=\tau_a \tau_b \tau_a(w) \tau_a \tau_b(f(a,c))\tau_a(f(b,ac))f(a,b\cdot ac)
\end{equation}
for all $w\in K$ and $a,b,c\in E$. Taking $w = 1$ gives \peqref{LBExtension1},
while taking $c = 1$ gives \peqref{LBExtension2}. Conversely, it is easy to
see that if both \peqref{LBExtension1} and \peqref{LBExtension2} hold, then
\peqref{LBExt-tmp} holds.

For (ii), we merely substitute $x=(u,a)$, $y=(v,b)$, and $z=(w,c)$ into
the associative law $x\cdot yz = xy\cdot z$. For (iii), then, it follows that
$Q$ is a group if and only if $E$ is a group and \peqref{RNExtension} holds
for all $a,b,c\in E$, $w\in K$. That \peqref{RNExtension}, universally
quantified, is equivalent to \peqref{GExtension1}
and \peqref{GExtension2} is proven similarly as in (i).

For part (iv), we plug $x=(u,a)$, $y=(v,b)$ into the commutative law $xy=yx$
to get that $(u,a)\in C(Q)$ if and only if $a\in C(E)$ and
\begin{equation}
\eqnlabel{InC-tmp}
v\tau_b(u) f(b,a) = u \tau_a(v) f(a,b)
\end{equation}
for all $v\in K$, $b\in E$. Taking $b = 1$ and rearranging gives \peqref{InC1},
while taking $v = 1$ and rearranging gives \peqref{InC2}. Conversely, it is easy to
see that if both \peqref{InC1} and \peqref{InC2} hold, then \peqref{InC-tmp} holds.
\end{proof}

In the next two sections, we will consider two special cases.
The extension $Q(K,E,\tau,f)$ is a \emph{semidirect product} if $f:E\times E\to
K$ satisfies $f(a,b)=1$ for every $a$, $b\in E$. In such a case, we denote the
resulting loop by $Q(K,E,\tau)$. Semidirect products of Bol loops were
considered in \cite{GR}.

When $A$, $B$ are loops, a map $\varphi:A\to B$ is a \emph{semihomomorphism} if
$\varphi(a\cdot ba) = \varphi(a)\cdot\varphi(b)\varphi(a)$ for every $a$, $b\in
A$.

\begin{corollary}
\corlabel{LBSemidirect}
Let $K$ be a group, $E$ a Bol loop, and $\tau:E\to\Aut(K)$ a map
satisfying $\tau_1=1$, and set $Q=Q(K,E,\tau)$. Then:
\begin{enumerate}
\item[(i)] $Q$ is a Bol loop if and only if $\tau$ is a semihomomorphism.
\item[(ii)] $(w,c)\in Q$ belongs to $\RNuc(Q)$ if and only if
$c\in \RNuc(E)$ and $\tau_{ab}(w) = \tau_a \tau_b(w)$ for all $a,b\in E$.
\item[(iii)] $Q$ is a group if and only if $E$ is a group and $\tau$ is a homomorphism.
\item[(iv)] $(u,a)\in Q$ belongs to $C(Q)$ if and only if $a\in C(E)$,
$\tau_a(v)=u^{-1}vu$ for every $v\in K$, and $u=\tau_b(u)$ for every $b\in E$.
\item[(v)] $C(Q)\subseteq \RNuc(Q)$ if and only if $C(E)\subseteq \RNuc(E)$.
\end{enumerate}
\end{corollary}

\begin{proof}
Parts (i), (ii), (iii), and (iv) follow immediately from specializing the corresponding
parts of Theorem \thmref{LB}. For (v), if $(w,c)\in C(Q)$, then by (iv), we trivially
have $\tau_{ab}(w) = \tau_a \tau_b(w)$ for all $a,b\in E$. If $C(E)\subseteq \RNuc(E)$,
then $c\in \RNuc(E)$, and so $(w,c)\in \RNuc(Q)$ by (ii). Conversely, if
$C(Q)\subseteq \RNuc(Q)$, then for $c\in C(E)$, $(1,c)\in \RNuc(Q)$, and so
$c\in \RNuc(E)$ by (ii).
\end{proof}

For $f\in\Aut(K)$, let $\Fix(f)=\{u\in K\,|\,f(u)=u\}$. For a
map $\tau:E\to\Aut(K)$, let $\Ker(\tau) = \setof{e\in E}{\tau_e=1}$ and
$\Fix(\tau) = \setof{u\in K}{u\in\Fix(\tau_e)\text{ for every }e\in E}$.

\begin{corollary}
\corlabel{KF}
Let $E$, $K$, $\tau$ and $Q=Q(K,E,\tau)$ be as in Corollary
$\corref{LBSemidirect}$. If both $E$ and $K$ are commutative, then
$C(Q) = \setof{(u,a)}{u\in\Fix(\tau)\text{ and }a\in\Ker(\tau)}$ and
$|C(Q)| = |\Fix(\tau)|\cdot |\Ker(\tau)|$.
\end{corollary}

The other special case of Theorem \thmref{LB} we consider is where the
``action'' $\tau : E\to \Aut(K)$ is trivial.
Denote by $\iota$ the map $\iota:E\to\Aut(K); a\mapsto 1$.

\begin{corollary}
\corlabel{LBTriv}
Let $K$ be a group, $E$ a Bol loop, $f:E\times E\to K$ a
cocycle, and set $Q=Q(K,E,\iota,f)$.
Then:
\begin{enumerate}
\item[(i)] $Q$ is a Bol loop if and only if $f(b,a)f(a,ba)\in Z(K)$
and
\begin{equation}
\eqnlabel{LBTriv}
f(b,a) f(a,ba)f(a\cdot ba,c) = f(a,c) f(b,ac) f(a,b\cdot ac)
\end{equation}
for every $a$, $b$, $c\in E$.
\item[(ii)] $(w,c)\in Q$ belongs to $\RNuc(Q)$ if and only if
$c\in \RNuc(E)$ and
\begin{equation}
\eqnlabel{RNTriv}
w\iv f(a,b) w f(ab,c) = f(b,c) f(a,bc)
\end{equation}
for all $a,b\in E$.
\item[(iii)] $Q$ is a group if and only if $E$ is a group, $f(a,b)\in Z(K)$,
and
\begin{equation}
f(b,c) f(a,bc) = f(a,b) f(ab,c)
\end{equation}
for every $a$, $b$, $c\in E$.
\item[(iv)] $(u,a)\in Q$ belongs to $C(Q)$ if and only if $a\in C(E)$,
$u\in Z(K)$, and $f(a,b)=f(b,a)$ for every $b\in E$.
\end{enumerate}
\end{corollary}

\begin{proof}
These claims follow immediately from specializing the corresponding
parts of Theorem \thmref{LB}.
\end{proof}

\section{Constructions based on semidirect products}
\seclabel{semi}

Moorhouse classified all nonassociative right Bol loops of order $16$, viz
\cite{Mo}. It turns out that among these $2049$ loops precisely $21$ have a
non-subloop commutant: $1$ of order $12$, and $20$ of order $16$. Among the
$20$ loops of order $16$, $19$ loops have commutant of order $6$, and $1$ loop
has commutant of order $4$.

We show in this subsection that precisely $3$ of the $21$ loops can be obtained
by a semidirect construction. All $21$ loops will be constructed in the next
section.

\begin{proposition}\prplabel{elem}
Let $K$ be a group, $E$ an elementary abelian $2$-group,
$\tau:E\to\mathrm{Aut}(K)$ a map such that $\tau_1=1$, $|\tau_e|=2$ for every
$e$, and $\langle \tau_E\rangle$ is a commutative subgroup of
$\mathrm{Aut}(K)$. Assume further that there are $a$, $b\in E$ such that
$\tau_a=\tau_b=1\ne\tau_{ab}$. Then:
\begin{enumerate}
\item[(i)] $\tau$ is a semihomomorphism but not a homomorphism,
\item[(ii)] $Q=Q(K,E,\tau)$ is a nonassociative Bol loop,
\item[(iii)] $C(Q)$ is not a subloop of $Q$,
\item[(iv)] $C(Q)\subseteq \RNuc(Q)$.
\end{enumerate}
\end{proposition}

\begin{proof}
We have $\tau_{a\cdot ba} = \tau_b$ since $E$ is an elementary abelian
$2$-group. On the other hand, $\tau_a\cdot \tau_b\tau_a = \tau_b$ since
$\langle\tau_E\rangle$ is commutative and $\tau_a$ is an involution. The
condition $\tau_a=\tau_b=1\ne\tau_{ab}$ guarantees that $\tau$ is not a
homomorphism. This proves (i). Then (ii) follows by Corollary \corref{LBSemidirect}.
Given $a$, $b\in E$ such that $\tau_a=\tau_b=1\ne\tau_{ab}$, note that
$(1,a)$, $(1,b)$ belong to $C(Q)$ but $(1,ab)$ does not, and so (iii) holds.
Finally, (iv) follows from Corollary \corref{LBSemidirect}(v) since $E$
is a group.
\end{proof}

\begin{example}\exmlabel{order12}
Let $E=\langle e_1,e_2\rangle$ be the elementary abelian $2$-group of order $4$,
and $K$ the cyclic group of order $3$, $\Aut(K)=\{1,\varphi\}$. Define
$\tau:E\to\Aut(K)$ by $\tau_1=\tau_{e_1}=\tau_{e_2}=1$,
$\tau_{e_1e_2}=\varphi$. Then $|\Ker(\tau)|=3$, $|\Fix(\tau)| =
|\Fix(\varphi)| = 1$. Hence $Q=Q(K,E,\tau)$ is a nonassociative
Bol loop of order $12$ with non-subloop commutant of order $3$.
\end{example}

\begin{example}\exmlabel{order16cyclic}
Let $E=\langle e_1,e_2\rangle$ be the elementary abelian $2$-group of order $4$,
and $K$ the cyclic group of order $4$, $\Aut(K)=\{1,\psi\}$. Define
$\tau:E\to\Aut(K)$ by $\tau_1=\tau_{e_1}=\tau_{e_2}=1$,
$\tau_{e_1e_2}=\psi$. Then $|\Ker(\tau)|=3$, $|\Fix(\tau)| = |\Fix(\psi)| = 2$.
Hence $Q=Q(K,E,\tau)$ is a nonassociative Bol
loop of order $16$ with non-subloop commutant of order $6$. It is easy to check
that $Q$ contains $9$ involutions.
\end{example}

\begin{example}\exmlabel{order16elem}
Assume that both $E=\langle e_1,e_2\rangle$ and $K=\langle k_1, k_2\rangle$ are
elementary abelian $2$-groups of order $4$. Define $\tau:E\to\Aut(K)$
by $\tau_1=\tau_{e_1}=\tau_{e_2}=1$, $\tau_{e_1e_2}: k_1\mapsto k_1$,
$k_2\mapsto k_1k_2$. Then $|\Ker(\tau)|=3$, $|\Fix(\tau)| =
|\Fix(\tau_{e_1e_2})| = 2$. Hence $Q=Q(K,E,\tau)$ is a nonassociative
Bol loop of order $16$ with non-subloop commutant of order $6$. It is easy
to check that $Q$ contains $13$ involutions.
\end{example}

\begin{lemma}
\lemlabel{E2Impossible}
Let $E$, $K$, $\tau$ be as in Corollary $\corref{LBSemidirect}$. If $|E|=2$ or
$|K|=2$ then $Q(K,E,\tau)$ is a group if and only if it is a Bol loop.
\end{lemma}

\begin{proof}
Let $E=\{1,e\}$, and assume that $\tau$ is a semihomomorphism. Then $\tau_e =
\tau_{eee} = \tau_e\tau_e\tau_e$ implies $\tau_{ee}=1=\tau_e\tau_e$, and hence
$\tau$ is a homomorphism.

If $|K|=2$ then $\Aut(K)=\{1\}$ and $\tau$ is a homomorphism.
\end{proof}

\begin{lemma}
Of the known Bol loops of order at most $16$ with non-subloop commutant, those constructed in Examples
$\exmref{order12}$, $\exmref{order16cyclic}$, $\exmref{order16elem}$ are the
only ones  obtained
by a nontrivial $(|E|>1$ and $|K|>1)$ application of the semidirect
construction of Corollary $\corref{LBSemidirect}$.
\end{lemma}

\begin{proof}
We rely on Moorhouse's classification \cite{Mo}; the caveat ``known'' in the
statement of the lemma is because the classification of the Bol loops of order
$16$ has not been independently verified. By Corollary \corref{order4k}, the only
possible orders less than or equal to $16$ for Bol loops with non-subloop commutants
are $8$, $12$, and $16$.
None of the Bol loops of order $8$ have
non-subloop commutant. (This also follows from Burn's classification of
Bol loops of order $8$ \cite{Burn1}.)

The loop of Example
\exmref{order12} is the only Bol loop of order $12$ with non-subloop
commutant, by the classification. (This also follows from Burn's
classification of Bol loops of order $4p$, $p$ an odd prime \cite{Burn2}.)

Assume that $Q=Q(K,E,\tau)$ is a Bol loop of order $16$ with non-subloop
commutant. By Lemma \lemref{E2Impossible}, we can assume that $|E|=4$ and
$|K|=4$. Let $k=|\mathrm{Ker}(\tau)|$, $f=|\mathrm{Fix}(\tau)|$. Since both $E$
and $K$ are abelian, $|C(Q)|=kf$ by Corollary \corref{KF}. By the
classification, the only possible values of $|C(Q)|$ are $4$ and $6$. If $k=4$
or $f=4$, $\tau_e=1$ for every $e\in E$ and hence $\tau$ is a homomorphism, a
contradiction.

If $K$ is cyclic, we have $f=2$ iff there is $e\in E$ such that $\tau_e$ is the
unique involution of $\mathrm{Aut}(\mathbb Z_4)$. If $K$ is elementary abelian,
we have $\mathrm{Aut}(K)\cong S_3$, and hence $f=2$ if and only if all
non-identity automorphism $\tau_e$ are equal to the same involution of
$\mathrm{Aut}(K)$.

Assume $|C(Q)|=4$. Then $k=f=2$. If $E=\langle e_1,e_2\rangle$ is elementary
abelian, we can assume that $\tau_1=\tau_{e_1}=1$ and $1\ne
\tau_{e_2}=\tau_{e_1e_2}$ is an involution. But then $\tau$ is a homomorphism,
a contradiction. If $E=\langle e\rangle$ is cyclic, then we can assume that
either $\tau_1=\tau_e=1$ and $1\ne \tau_{e^2}=\tau_{e^3}$ is an involution,
which results in $1\ne
\tau_{e^3}\tau_{e^2}\tau_{e^3}=\tau_{e^3e^2e^3}=\tau_1=1$; or we can assume
that $\tau_1=\tau_{e^2}=1$ and $1\ne \tau_e=\tau_{e^3}$ is an involution, which
means that $\tau$ is a homomorphism.

Now assume that $|C(Q)|=6$. Then $k=3$, $f=2$. If $E=\langle e_1,e_2\rangle$ is
elementary abelian, we can assume that $1=\tau_{e_1}=\tau_{e_2}$ and $1\ne
\tau_{e_1e_2}$ is an involution. Since there is a unique involution in
$\mathrm{Aut}(\mathbb Z_4)$ and since $\mathrm{Aut(Aut(K))}\cong S_3$ acts
transitively on the involutions of $\mathrm{Aut(K)}\cong S_3$ when $K$ is
elementary abelian, this case yields the loops obtained in Examples
\exmref{order16cyclic} and \exmref{order16elem}. Finally assume that $E=\langle
e\rangle$ is cyclic. Then we can assume that either $\tau_e=\tau_{e^2}=1$ and
$1\ne \tau_{e^3}$ is an involution, which yields $1\ne
\tau_{e^3}=\tau_{eee}=\tau_e\tau_e\tau_e = 1$; or we can assume that $\tau_e =
\tau_{e^3} = 1\ne\tau_{e^2}$, which yields
$1\ne\tau_{e^2} = \tau_e\tau_{e^2}\tau_e = \tau_{ee^2e} = 1$, a contradiction.
\end{proof}

Note that if the commutant $C(Q)$ of a Bol loop $Q$ has order $2$, say,
$C(Q) = \{1,a\}$, then $C(Q)$ is a subloop. By contrast, we have the following.

\begin{proposition}
\prplabel{orderc}
For each $k>2$, there exists a Bol loop with non-subloop commutant of order $k$.
\end{proposition}

\begin{proof}
Pick $n$ such that $2^n>k$. Let $E$ be the elementary abelian $2$-group of
order $2^n$, and let $K$ be the cyclic group of order $3$, thus
$\mathrm{Aut}(K)=\{1,\varphi\}$. For some $a\ne 1\ne b\ne a$ in $E$, let
$\tau_1=\tau_{a}=\tau_{b}=1$, $\tau_{ab}=\varphi$. Choose the remaining $2^n-4$
automorphisms $\tau_e$ of $K$ arbitrarily, but in such a way that
$|\mathrm{Ker}(\tau)|=k$. Then $Q=Q(K,E,\tau)$ is a nonassociative Bol
loop by Proposition \prpref{elem}. Moreover, since
$|\mathrm{Fix}(\tau)|=|\mathrm{Fix}(\varphi)|=1$ and both $E$ and $K$ are abelian,
$|C(Q)| = k$ by Corollary \corref{KF}.
\end{proof}

\begin{proposition}
\prplabel{order4n}
For each $n>2$, there exists a Bol loop of order $4n$ with non-subloop commutant.
\end{proposition}

\begin{proof}
Let $E=\langle e_1,e_2\rangle$ be the elementary abelian group of order $4$, and
let $K$ be the cyclic group of order $n$. Then $\psi:k\mapsto k^{-1}$ is a
non-identity involutory automorphism of $K$. Set
$\tau_1=\tau_{e_1}=\tau_{e_2}=1$, $\tau_{e_1e_2}=\psi$. By Proposition
\prpref{elem}, $Q=Q(K,E,\tau)$ is a nonassociative Bol loop of order $4n$
with non-subloop commutant.
\end{proof}

\section{Constructions based on extensions}
\seclabel{trivact}

In this section, we will use additive notation for abelian groups.
As an immediate consequence of Corollary \corref{LBTriv} we get:

\begin{lemma}
\lemlabel{fabfba}
Let $K$ be an abelian group, $E$ be a group, $f:E\times E\to K$ a cocycle, and
$Q=Q(K,E,\iota,f)$. Then:
\begin{enumerate}
\item[(i)] $(w,c)\in Q$ belongs to $\RNuc(Q)$ if and only if
$f(a,b) + f(ab,c) = f(b,c) + f(a,bc)$ for all $a,b\in E$,
\item[(ii)] $(u,a)\in Q$ belongs to $C(Q)$ if and only if
$f(a,b)=f(b,a)$ for every $b\in E$.
\end{enumerate}
\end{lemma}

\begin{lemma}\lemlabel{E2K2}
Let $E$ and $K$ be elementary abelian $2$-groups, $f:E\times E\to K$ a cocycle,
and $Q=Q(K,E,\iota,f)$. Then:
\begin{enumerate}
\item[(i)] $Q$ is a Bol loop if and only if
\begin{eqnarray}
\eqnlabel{triv-exp2-1}
f(a, a + c) &=& f(a,a) + f(a,c), \\
\eqnlabel{triv-exp2-2}
f(a, b + c) + f(a,b) + f(a,c) &=& f(b, a + c) + f(b,a) + f(b,c)
\end{eqnarray}
for all $a,b,c\in E$.
\item[(ii)] If there exist $a,b,c\in E$ such that $f(a+b,c)\ne f(c,a+b)$ and
$f(a,d)=f(d,a)$, $f(b,d)=f(d,b)$ for every $d\in E$, then $C(Q)$ is not a
subloop of $Q$.
\end{enumerate}
\end{lemma}

\begin{proof}
We freely use that $E$ and $K$ are of exponent $2$.
In additive notation, the cocycle identity \peqref{LBTriv} is
\begin{equation}
\eqnlabel{triv-exp2-tmp}
    f(b,a) + f(a,b+a) + f(b,c) = f(a,c) + f(b,a+c) + f(a,b+a+c) .
\end{equation}
Taking $b = a$, we get \peqref{triv-exp2-1}, and applying
\peqref{triv-exp2-1} to \peqref{triv-exp2-tmp}, we get \peqref{triv-exp2-2}.
Conversely, it is easy to see that \peqref{triv-exp2-1} and \peqref{triv-exp2-2}
imply \peqref{triv-exp2-tmp}. This establishes (i).

Assume that $a$, $b$, $c$ are as in (ii). Then $(0,a)$, $(0,b)$ belong to
$C(Q)$ by Lemma \lemref{fabfba}. By the same Lemma, $(0,a)(0,b)=(f(a,b),a+b)$
does not belong to $C(Q)$.
\end{proof}

In case $E$ is an abelian group,
we say that the cocycle $f:E\times E\to K$ is \emph{right additive} if
$f(a,b+c)=f(a,b)+f(a,c)$ for every $a,b,c\in E$.

\begin{lemma}\lemlabel{E2K2-b}
Let $E$ and $K$ be elementary abelian $2$-groups, $f:E\times E\to K$ a
right additive cocycle, and $Q=Q(K,E,\iota,f)$. Then $Q$ is a Bol loop, and
\begin{enumerate}
\item[(i)] $(w,c)\in Q$ belongs to $\RNuc(Q)$ if and only if the mapping
$E\to K; a\mapsto f(a,c)$ is a homomorphism,
\item[(ii)] $C(Q)\subseteq \RNuc(Q)$.
\end{enumerate}
\end{lemma}

\begin{proof}
That $Q$ is a Bol loop follows immediately from Lemma \lemref{E2K2}(ii) and
right additivity. Again using right additivity,
Lemma \lemref{fabfba}(i) reduces to $(w,c)\in \RNuc(Q)$ if and only if
$f(a + b,c) = f(a,c) + f(b,c)$ for all $a,b\in E$, and this establishes (i).
Finally, if $(w,c)\in C(Q)$, then by Lemma \lemref{fabfba}(ii) and
right additivity, $f(a+b,c) = f(c,a+b)
= f(c,a) + f(c,b) = f(a,c) + f(b,c)$ for all $a,b\in E$, and so
$(w,c)\in \RNuc(Q)$ by (i).
\end{proof}

When $K=\{0,1\}\cong \mathbb{Z}_2$ and $E$ is an elementary
abelian $2$-group, then $E$ is a vector space over $K$ and
a cocycle $f:E\times E\to K$ is a form satisfying
$f(0,a) = f(a,0) = 0$. As usual, we say
that $g:E\times E\to K$ is \emph{equivalent} to $f$ if there is
$\varphi\in\mathrm{Aut}(E)$ such that $f(a,b) = g(\varphi(a),\varphi(b))$ for
every $a$, $b\in E$.

\begin{lemma} Let $E$ be a vector space over $K=\{0,1\}$ with basis $B=\{e_1$,
$\dots$, $e_n\}$. Let $c:E\times B\to K$ be a map satisfying $c(0,e_i)=0$ for
every $1\le i\le n$. Then there is a unique right additive cocycle $f:E\times
E\to K$ such that $f(e,e_i)=c(e,e_i)$ for every $e\in E$, $e_i\in B$.
\end{lemma}
\begin{proof} The map $f(a,\underline{\phantom{a}}\,):E\to K$,
$a\mapsto f(a,e)$ is a homomorphism for every $a$ if and only if $f$ is right
additive.
\end{proof}

In the situation of the previous Lemma, we say that $f$ is \emph{associated}
with $c$.

\begin{proposition}\prplabel{f}
Let $E$ be a vector space over $K=\{0,1\}$ with basis $B=\{e_1$, $\dots$,
$e_n\}$. Let $c:E\times B\to K$ be a map satisfying $c(0,e_i)=0$, $c(e_1,e_i) =
c(e_i,e_1)$, $c(e_2,e_i)=c(e_i,e_2)$ for every $1\le i\le n$, and
$c(e_1+e_2,e_3)\ne c(e_1,e_3)+c(e_2,e_3)$. Assume furthermore that $c$ is such
that the right additive cocycle $f_c=f:E\times E\to K$ associated with $c$
satisfies $f(e_1,e)=f(e,e_1)$, $f(e_2,e)=f(e,e_2)$ for every $e\in E$. Then
$Q=Q(K,E,\iota,f)$ is a Bol loop with non-subloop commutant.

Moreover, if $g:E\times E\to K$ is a right additive cocycle satisfying the
assumptions of Lemma \emph{\lemref{E2K2}(ii)}, then $g$ is equivalent to $f_c$ with some
choice of $c$ as above.
\end{proposition}

\begin{proof}
The conditions $c(e_1+e_2,e_3)\ne c(e_1,e_3)+c(e_2,e_3)$, $f(e,e_1)=f(e_1,e)$,
and $f(e,e_2)=f(e_2,e)$ guarantee that $f(e_1+e_2,e_3)\ne f(e_1,e_3)+f(e_2,e_3)
= f(e_3,e_1)+f(e_3,e_2)=f(e_3,e_1+e_2)$. Hence $f$ satisfies the assumptions of
Lemma \lemref{E2K2}(ii), and $Q=Q(K,E,\iota,f)$ is a Bol loop with
non-subloop commutant.

Let $g:E\times E\to K$ be a right additive cocycle with $a$, $b$, $c\in E$ such
that $g(a,d)=g(d,a)$, $g(b,d)=g(d,b)$ for every $d\in E$, and such that
$g(a+b,c)\ne g(c,a+b)$. Then $g(a+b,c)\ne g(c,a)+g(c,b)=g(a,c)+g(b,c)$.

It is clear that neither $a$ nor $b$ nor $c$ can be equal to $0$, and that
$a\ne b$. Moreover, $c\ne a$ (else $g(a+b,c) = g(a+b,a) = g(a,a+b) = g(a,a) +
g(a,b) = g(a,a) + g(b,a) = g(a,c)+g(b,c)$), $c\ne b$ (by a similar argument),
and $c\ne a+b$ (else $g(a+b,c)=g(c,a+b)$). Hence $a$, $b$, $c$ are linearly
independent, and there is an automorphism of $E$ that maps $(a,b,c)$ to
$(e_1,e_2,e_3)$.
\end{proof}

\begin{corollary}\corlabel{f}
Let $E$, $K$, $B$, $n$, $c$ and $f$ be as in Proposition \emph{\prpref{f}}. Then we
are free to choose $(2^n-4)(n-2)+3n-4$ of the values of $c$.
\end{corollary}
\begin{proof} We can choose $c(e_1,e_i)$ for every $1\le i\le n$. Then
$c(e,e_1)$ is determined for every $e\in E$ by the condition
$f(e,e_1)=f(e_1,e)$. We can then choose $c(e_2,e_i)$ for every $2\le i\le n$,
hence determining $c(e,e_2)$ for every $e$ by the condition
$f(e,e_2)=f(e_2,e)$. Since $c(e_1+e_2,e_3)\ne c(e_1,e_3)+c(e_2,e_3)$, the value
$c(e_1+e_2,e_3)$ is determined. But we can choose $c(e_1+e_2,e_i)$ for every
$4\le i\le n$. Finally for $e$ not in the subspace $\langle e_1,e_2\rangle$, we
are free to choose $c(e,e_i)$ for every $3\le i\le n$.
\end{proof}

\begin{example}\exmlabel{c1c9}
Let $E=\langle e_1,e_2,e_3\rangle$ be a $3$-dimensional vector space over
$K=\{0,1\}$. According to Corollary \corref{f}, we are free to choose
$(2^3-1)(3-2)+3\cdot 3 - 4 = 9$ values of $c:E\times \{e_1,e_2,e_3\}\to K$ in
order to uniquely determine the associated right additive cocycle $f:E\times
E\to K$ such that $Q=Q(K,E,\iota,f)$ is a Bol loop with non-subloop
commutant.

These nine choices are as follows:
\begin{displaymath}
    \begin{array}{r|ccc}
         c& e_1& e_2& e_3\\
         \hline
        e_1&c_1&c_2&c_4\\
        e_2&.&c_3&c_5\\
        e_1+e_2&.&.&.\\
        e_3&.&.&c_6\\
        e_1+e_3&.&.&c_7\\
        e_2+e_3&.&.&c_8\\
        e_1+e_2+e_3&.&.&c_9\\
    \end{array}
\end{displaymath}
The resulting loop of order $16$ will be denoted by
$Q(c_1,c_2,c_3,c_4,c_5,c_6,c_7,c_8,c_9)$.
\end{example}

\begin{proposition}\prplabel{all16}
All Bol loops of order $16$ with non-subloop
commutant are isomorphic to the loop $Q(c_1,\cdots,c_9)$ with some choice of
$c_1$, $\dots$, $c_9\in K=\{0,1\}$, except for one loop.
\end{proposition}
\begin{proof}
We have verified by computer, using the package LOOPS \cite{LOOPS}, that the
following $19$ Bol loops are pairwise non-isomorphic:\
\begin{small}
\begin{align*}
   &Q(0,0,0,0,0,0,0,0,0),\quad Q(0,0,0,0,0,0,0,0,1),\quad Q(0,0,0,0,0,0,0,1,1),\\
   &Q(0,0,0,0,0,0,1,1,0),\quad Q(0,0,0,0,0,0,1,1,1),\quad Q(0,0,0,0,0,1,1,1,1),\\
   &Q(0,0,1,0,0,0,0,0,0),\quad Q(0,0,1,0,0,0,0,0,1),\quad Q(0,0,1,0,0,0,0,1,1),\\
   &Q(0,0,1,0,0,0,1,0,0),\quad Q(0,0,1,0,0,0,1,0,1),\quad Q(0,0,1,0,0,0,1,1,0),\\
   &Q(0,0,1,0,0,1,1,0,0),\quad Q(1,0,1,0,0,0,0,0,0),\quad Q(1,0,1,0,0,0,0,0,1),\\
   &Q(1,0,1,0,0,0,0,1,0),\quad Q(1,0,1,0,0,0,0,1,1),\quad Q(1,0,1,0,0,0,1,1,0),\\
   &Q(1,0,1,0,0,1,0,0,1).
\end{align*}
\end{small}
There is only one additional Bol loop with non-subloop commutant,
according to Moorhouse's classification.
\end{proof}

Here is the unique Bol loop $Q$ of order $16$ with non-subloop commutant
not obtained in Proposition \prpref{all16}:
\begin{small}
\begin{displaymath}
\begin{array}{r|cccc|cccc|cccc|cccc}
&K&&&&Ke_1&&&&Ke_2&&&&Ke_1e_2&&&\\
\hline
1&1&2&3&4&5&6&7&8&9&10&11&12&13&14&15&16\\
k_1&2&1&4&3&6&5&8&7&10&9&12&11&14&13&16&15\\
k_2&3&4&1&2&7&8&5&6&12&11&10&9&16&15&14&13\\
k_1k_2&4&3&2&1&8&7&6&5&11&12&9&10&15&16&13&14\\
\hline
Ke_1&5&6&7&8&1&2&3&4&13&14&15&16&9&10&11&12\\
&6&5&8&7&2&1&4&3&16&15&14&13&12&11&10&9\\
&7&8&5&6&3&4&1&2&15&16&13&14&11&12&9&10\\
&8&7&6&5&4&3&2&1&14&13&16&15&10&9&12&11\\
\hline
Ke_2&9&10&11&12&13&14&15&16&1&2&3&4&5&6&7&8\\
&10&9&12&11&14&13&16&15&2&1&4&3&6&5&8&7\\
&11&12&9&10&15&16&13&14&3&4&1&2&7&8&5&6\\
&12&11&10&9&16&15&14&13&4&3&2&1&8&7&6&5\\
\hline
Ke_1e_2&13&14&15&16&9&10&11&12&5&6&7&8&1&2&3&4\\
&14&13&16&15&10&9&12&11&6&5&8&7&2&1&4&3\\
&15&16&13&14&11&12&9&10&7&8&5&6&3&4&1&2\\
&16&15&14&13&12&11&10&9&8&7&6&5&4&3&2&1
\end{array}
\end{displaymath}
\end{small}
Note that $\LNuc(Q) = Z(Q) = \{1\}$, $\RNuc(Q) = \{1, 2, 3 ,4, 5, 6, 7, 8\}$ is
an elementary abelian $2$-group, and $C(Q) = \{1,2,5,7\}$. Also note that
$\RNuc(Q) = \sbl{C(Q)}$.

The fact that $N_\lambda(Q)$ is trivial implies that $Q$ cannot be obtained by
any extension $(\eqnref{AllExts})$ of Theorem \thmref{AllExts}. In
\cite{KiechleNagy}, Kiechle and Nagy developed a theory of extensions of
involutory Bol loops, and constructed all involutory Bol loops of
order $16$, three of which happen to have a trivial center. Since our loop $Q$
is involutory and has trivial center, it is one of the three loops mentioned in
\cite[Corollary 7]{KiechleNagy}.

We conclude this section with an explicit construction of $Q$. Let $K=\langle
k_1, k_2\rangle$ and $E=\langle e_1,e_2\rangle$ be elementary abelian
$2$-groups of order $4$. For every $(a,b)\in E\times E$ we define an
automorphism $\psi_{a,b}$ of $K$. Namely: $\psi_{1,e_2}=\psi_{1,e_1e_2}$
satisfies $k_1\mapsto k_1$, $k_2\mapsto k_1k_2$,
$\psi_{e_1,e_1}=\psi_{e_1,e_1e_2}$ satisfies $k_1\mapsto k_1k_2$, $k_2\mapsto
k_2$, and all other automorphisms $\psi_{a,b}$ are trivial. Then $Q$ is
isomorphic to $K\times E$ with multiplication $(u,a)(v,b) =
(\psi_{a,b}(u)v,ab)$, as is easily seen from the multiplication table of $Q$.
(For the convenience of the reader, we have subdivided the multiplication table
of $Q$ into subsquares corresponding to the cosets of $K$, labeled the cosets
of $K$, and also labeled the elements in one of the cosets.)

\section{Acknowledgement}

Our investigations were aided by the automated reasoning tools OTTER
\cite{Mc-Otter} and Prover9 \cite{Mc-Prover9}, and by the finite model builder
Mace4 \cite{Mc-Mace}.

\end{document}